\documentclass[11pt,a4paper]{amsart}

\usepackage{amsfonts}
\usepackage{amssymb}
\usepackage{mathrsfs} 

\theoremstyle{plain}
\newtheorem{theorem}{Theorem}

\newtheorem{corollary}[theorem]{Corollary}

\theoremstyle{definition}
\newtheorem{definition}[theorem]{Definition}

\newtheorem{remark}[theorem]{Remark}
\newtheorem{remarks}[theorem]{Remarks}

\newcommand{\dA}{\mathrm{det}A}

\newcommand{\Order}{\mathcal{O}}
\newcommand{\vol}{\mathrm{Vol}}
\newcommand{\area}{\mathrm{Area}}
\newcommand{\length}{\textrm{Length}}
\newcommand{\llangle}{\left\langle}
\newcommand{\rrangle}{\right\rangle}

\newcommand{\dd}{\mathrm{d}}
\newcommand{\I}{\mathrm{I}}
\newcommand{\II}{\mathrm{I\!I}}

\newcommand{\Hess}{\mathrm{Hess}}
\newcommand{\F}{\mathcal{F}}
\newcommand{\M}{\mathcal{M}}
\newcommand{\E}{\mathcal{E}}
\newcommand{\deltatnought}{\left.\frac{\partial}{\partial t}\right\vert_{t=0}}
\newcommand{\grad}{\textrm{grad}}
\newcommand{\DD}{\overline{\mathrm{D}}} 

\title[On the Area Functional of the Second Fundamental Form]{On the Area Functional\\ of the Second Fundamental Form\\ of  Ovaloids.}

\author[S. Verpoort]{Steven Verpoort${}^{1}$\\ Katholieke Universiteit Leuven, Belgium.}

\email{steven.verpoort@wis.kuleuven.be}

\begin{document}

\maketitle

\begin{abstract}
The expression for the variation of the area
functional of the second fundamental form of a hypersurface in a
Euclidean space involves the so-called ``\textit{mean curvature of the second fundamental form}.'' Several new characteristic properties of (hyper)spheres, in which the mean curvature of the second fundamental form occurs, are given. In particular, it is shown that the spheres are the only ovaloids which are a critical point of the area functional of the second fundamental form under various constraints.
\end{abstract}


\section{Introduction.}

The intention of the article at hand is to contribute to the theory of hypersurfaces in a Euclidean space,
of which the second fundamental form is positive-definite and accordingly can be seen as an abstract Riemannian metric.
\hspace{50cm}\footnote{The author was partially supported by the Research
Foundation -- Flanders (project G.0432.07).}

\subsection{Some Results on the Geometry of the Second Fundamental Form.}

The \textit{intrinsic} geometry of this abstract metric and its influence on the shape of the original hypersurface
have already been the object of extensive study. For example, it has been proved by R.\ Schneider that the
hyperspheres are the only such compact hypersurfaces of which the second fundamental form has
constant Riemannian curvature  \cite{schneider1972}.
In extension, it has been demonstrated that only hyperspheres can
satisfy certain further relations between the sectional curvature of the second fundamental form and the curvatures of the original compact hypersurface (e.g., \cite{baik1987,hasanis1980,koufo1977,koutroufiotis1974,simon,stamou1978}).

The notion of mean curvature, which belongs to the \textit{extrinsic} geometry of the hypersurface, and which can be characterised as a measure for the rate of area growth under deformations of the hypersurface, can be tailored to the geometry of the second fundamental form: a function which measures the rate at which the total area of the hypersurface, as surveyed in the geometry of its second fundamental form, changes under a deformation is called \textit{the mean curvature of the second fundamental form}. This concept has been introduced by E.\ Gl\"assner for surfaces in $\mathbb{E}^3$ \cite{glassner,glassnersimon} (see also, e.g., \cite{manhart1983} for hypersurfaces in $\mathbb{E}^{m+1}$).

\subsection{Overview.}

The notation which will be adopted is outlined in \S\,\ref{sec:not} below.
In \S\,\ref{sec:IIvar}, the expression for the variation of the area of the second fundamental form is given and some useful remarks regarding the mean curvature of the second fundamental form will be made.

As will be mentioned in the beginning of \S\,\ref{sec:characHii}, several characterisations of Euclidean hyperspheres in which the mean curvature of the second fundamental form occurs have been found already. Two similar new results are contained in this section.

In \S\,\ref{sec:int} an integral formula is derived which will be helpful in the next 
\S\,\ref{sec:varcharac}. Here it is shown that the spheres are the only ovaloids which are a critical point of the area functional of the second fundamental form under various constraints (Theorems~\ref{thm:firstvariationalproblem}, \ref{thm:secondvariationalproblem} and \ref{thm:thirdvariationalproblem}).


\section{Notation.}
\label{sec:not}

For a manifold $M$, the collection of all vector fields on $M$ will be denoted by $\mathfrak{X}(M)$ whereas $\mathfrak{F}(M)$ stands for the collection of all real-valued functions on $M$. 

If $f:M\rightarrow\mathbb{R}$ is a real-valued function on a Riemannian manifold $(M,g)$ with Levi-Civita connection $\nabla$, the \textit{Hessian operator} of $f$ is defined as
\[
\mathrm{Hs}_f : \mathfrak{X}(M)\rightarrow \mathfrak{X}(M) : V \mapsto \nabla_V \left(\grad f \right)\,.
\]
The \textit{Hessian} of $f$ (notation $\mathrm{Hess}_f$) is the  $(0,2)$-tensor which is metrically equivalent to $\mathrm{Hs}_f$. The \textit{Laplacian} of $f$ is defined as $\Delta f=\textrm{trace}\left(\mathrm{Hs}_f\right)$.

For a hypersurface $M\subseteq\mathbb{E}^{m+1}$,
the \textit{first fundamental form} will be denoted by $\I$ or $g$, and is the restriction of the Euclidean scalar product $\llangle\cdot\,,\,\cdot\rrangle$
to the tangent spaces of $M$.
The \textit{shape operator} $A$ is defined with respect to a unit normal vector field $N$ by
\[
A: \mathfrak{X}(M) \rightarrow \mathfrak{X}(M):V\mapsto A(V)=-\DD_V N\,,
\]
where $\DD$ is the
standard connection on $\mathbb{E}^{m+1}$. The \textit{second fundamental form} is given by $\II(V,W)=\langle A(V),W\rangle$ for
$V,W\in\mathfrak{X}(M)$, and the \textit{mean curvature} is equal to $H=\frac{1}{m}\textrm{trace}\,A$. The \textit{Gaussian curvature} of a surface in $\mathbb{E}^3$ is denoted by $K$.

A hypersurface $M\subseteq\mathbb{E}^{m+1}$ will be called
\textit{locally strongly convex} if the second fundamental form $\II$ is positive-definite.
A \textit{hyperovaloid} (resp.\ \textit{ovaloid}) is a compact, locally strongly convex hypersurface $M\subseteq\mathbb{E}^{m+1}$ (resp.\ $\mathbb{E}^{3}$).

The second fundamental form furnishes an abstract Riemannian metric on a locally strongly convex hypersurface $M$, and an index $\II$ will indicate that a geometrical object is defined with respect to $(M,\II)$.
It should be noted that we need to choose the \textit{interior} unit normal vector field on a hyperovaloid in order to make the second fundamental form positive-definite.


\section{The Variation of the Area of the Second Fundamental Form.}
\label{sec:IIvar}

The letter $\E$ will designate the set
of all locally strongly convex hypersurfaces in $\mathbb{E}^{m+1}$.
Our first concern is the infinitesimal behaviour of the area functional associated to the second fundamental form:
\[
\area_{\II} : \E \rightarrow \mathbb{R} : M\mapsto
\area_{\II}(M)=\int_M\,\dd\Omega_{\II} = \int_M \sqrt{\dA}
\,\dd\Omega\,.
\]

By a \textit{deformation} of a hypersurface
$M\subseteq\mathbb{E}^{m+1}$ will be understood a smooth mapping (for some
$\varepsilon > 0$)
\[
\mu: \,\left]\,-\varepsilon\, ,\, \varepsilon \,\right[ \,\times M \rightarrow
\mathbb{E}^{m+1} : (t,n) \mapsto \mu_t(n)\,,
\]
such that for some compact set $\M\subseteq M$
\[
\left\{
\begin{array}{l}
\textrm{for all $n\in M$,\quad} \mu_0(n)=n \,;\\
\textrm{for all $n\in M\setminus\M$ and all $t\in\,\left]-\varepsilon \,,\, \varepsilon\,\right[$\,,\quad} \mu_t(n)=n \,.\\
\end{array}
\right.
\]
The $\mathbb{E}^{m+1}$-valued vector field
sending $n\in M$ to $\displaystyle\deltatnought\!\mu_t(n)$
will be called the \textit{deformation vector field}.

The proof of the following Theorem, which was given by E.\ Gl\"assner in \cite{glassner} in dimension $m=2$ and in \cite{verpoort2008} for $m\geqslant 2$, has been omitted here.

\begin{theorem}[E.\ Gl\"assner]
\label{deltaAII} 
Let $M$ be a locally strongly convex hypersurface in $\mathbb{E}^{m+1}$.
The variation of the area functional of the second fundamental
form along a deformation $\mu$ with deformation
vector field $X$ is given by
\[
\deltatnought {\area}_{\II}(\mu_t(M))  =
\int_{M} \left\{
\frac{m}{2}H
+\frac{1}{4}\Delta_{\II}(\log \dA)\right\}\langle -N,X\rangle \,\dd\Omega_{\II}\,.
\]
\end{theorem}

\begin{definition}
Let $M$ be a locally strongly convex hypersurface in $\mathbb{E}^{m+1}$. The \textit{mean
curvature of the second fundamental form} $H_{\II}$ is defined
by
\[
H_{\II} = \frac{m}{2}H
+\frac{1}{4}\Delta_{\II}(\log \dA)\,.
\]
If $H_{\II}=0$, the surface will be
called \textit{$\II$-minimal}.
\end{definition}

\begin{remark}
There do not exist $\II$-minimal hyperovaloids in $\mathbb{E}^{m+1}$.
\end{remark}

\begin{remarks}
This terminology was introduced by E. Gl\"assner \cite{glassner,glassnersimon} for surfaces in $\mathbb{E}^3$ (see also, e.g., \cite{manhart1983} for hypersurfaces in $\mathbb{E}^{m+1}$).
There is some inaccuracy in the terminology, since a critical point of ${\area}_{\II}$ is not necessarily a minimum: for $m=2$, M.\ Wiehe \cite[remark 4.4]{wiehe1998} has shown that the second-order variation of $\area_{\II}$ is \emph{negative-definite} for locally strongly convex surfaces in $\mathbb{E}^3$ with $H_{\II}=0$.
Despite what might be suggested by the name \textit{mean curvature of the
second fundamental form}, the quantity $H_{\II}$ is not determined by the second fundamental form $\II$ solely.
This is, of course, in contrast with the \textit{scalar curvature of the second fundamental form}, which will be denoted by $S_{\II}$
and is given by \cite{schneider1972}
\begin{equation}
\label{eq:SII}
S_{\II} = m(m-1) H + \mathcal{P} - \frac{\II(\grad_{\II} \textrm{det}A,\grad_{\II} \textrm{det}A)}{4(\textrm{det}A)^2}
\end{equation}
for a non-negative function $\mathcal{P}$. If $M$ is an ovaloid in three-dimensional Euclidean space,
we will denote the Gaussian curvature of $(M,\II)$ as $K_{\II}=\frac{1}{2}S_{\II}$.
\end{remarks}

\begin{remark}
F.\ Manhart \cite{manhart1983} has introduced the so-called $\II$-normal field 
\[
N_{\II} = \sqrt{\dA}\, N - \grad_{\II}\sqrt{\dA}
\]
along a hyperovaloid $M\subseteq\mathbb{E}^{m+1}$. 
Adopting the  language of \cite{simon_schwenk_viesel_1991}, this is a relative normal vector field with relative shape operator
\[
A_{(N_{\II})} : \mathfrak{X}(M) \rightarrow \mathfrak{X}(M) : V \mapsto 
A_{(N_{\II})} (V) = - \DD_V N_{\II}\,,
\]
the trace of which is given by $2\,\sqrt{\dA}\,H_{\II}.$
\end{remark}

\section{Characterisations of Euclidean Hyperspheres.}
\label{sec:characHii}

Several characterisations of the Euclidean hyperspheres in which the function $H_{\II}$ occurs have been found already.
For example, it is known that the hyperspheres are 
the only hyperovaloids in $\mathbb{E}^{m+1}$ (with $m\geqslant 2$) which satisfy $2(m-1)H_{\II} \geqslant S_{\II}$ 
(see \cite[Satz 3.7, with $\alpha=\frac{1}{2}$]{manhart1989a}). Furthermore, an ovaloid in $\mathbb{E}^3$ of which it is known that any of the functions
$K_{\II}-\sqrt{K}$, $H_{\II}-H$, or $H_{\II}-K_{\II}$  does not change sign, is a sphere (see \cite{koutroufiotis1974}, \cite{stamou1981} and \cite[Korollar 1, $\{$(b), with $\alpha=0\}$ along with $\{$(c), with $\alpha=\frac{1}{2}\}$ ]{stamou2003}).

The two following theorems, the first one of which generalises results of \cite{stamou1981}, give similar characterisations of hyperspheres. 
\begin{theorem}
\label{thm:HIIenH}
If a hyperovaloid $M\subseteq\mathbb{E}^{m+1}$ satisfies either $H_{\II}=\frac{m}{2}H+f(\dA)$
or $H_{\II}=H f(\dA)$ for an increasing function $f:\mathbb{R}\rightarrow\mathbb{R}$, then $M$ is a hypersphere.
\end{theorem}
\begin{proof}[of Theorem~\ref{thm:HIIenH}]
 Let $p_{+}$ and $p_{-}$ be two points where $\dA$ attains its maximum and minimum, respectively.
In the first case, we necessarily have
\begin{eqnarray*}
\left.4 f(\dA)\right\vert_{p_{+}} &=& \left.\Delta_{\II}(\log\dA)\right\vert_{p_{+}} \leqslant
0 \leqslant \left.\Delta_{\II}(\log\dA)\right\vert_{p_{-}} \\
&=&\left.4 f(\dA)\right\vert_{p_{-}} \leqslant \left.4 f(\dA)\right\vert_{p_{+}}\,.
\end{eqnarray*}
This shows that $\dA$ is constant, which is only possible if $M$ is a hypersphere \cite{suss1929a}. 

In the second case,
\[
\left. \big( H f(\dA) \big)\right\vert_{p_{+}}=H_{\II}(p_{+})=\left.\left(\frac{m}{2}H + \frac{1}{4}\Delta_{\II}(\log\dA)\right)\right\vert_{p_{+}}\leqslant \frac{m}{2}H(p_{+})
\]
and consequently, every  $p\in M$ satisfies 
\[
\left.f(\dA)\right\vert_{p}\leqslant\left.f(\dA)\right\vert_{p_{+}}\leqslant \frac{m}{2}\,.
\] 
But this means $\Delta_{\II} (\log\dA)\leqslant 0$, and $M$ is a hypersphere.
\end{proof}

\begin{theorem}
\label{thm:ineq}
If $M\subseteq\mathbb{E}^{m+1}$ is a hyperovaloid, then there holds
\[
\int H_{\II}\,\dd\Omega \geqslant \frac{m}{2}\int H \dd\Omega
\]
with equality if and only if $M$ is a hypersphere.
\end{theorem}
\begin{proof}[of Theorem~\ref{thm:ineq}]
This is an immediate consequence of the relation
\[
\Delta_{\II} (\log \dA) = \mathrm{div}\left(\grad_{\II}(\log \dA)\right)+\frac{1}{2}\II\left(\grad_{\II}(\log \dA)\rule{0pt}{12pt},\grad_{\II}(\log \dA)\right)\,,
\]
which can easily be obtained by direct computation.
\end{proof}


\section{An Integral Formula.}
\label{sec:int}

It is possible to adapt an integral formula of H.\ Minkowski to the geometry of the second fundamental form. By $P$ shall be denoted the position vector field of $\mathbb{E}^{m+1}$ with respect to an arbitrary origin, and $\rho$ will stand for the support function
\[
\rho \equiv \llangle -N, P \rrangle : M \rightarrow \mathbb{R}\,.
\]

\begin{theorem}
\label{thm:integral}
The following integral formula holds on any hyperovaloid $M\subseteq\mathbb{E}^{m+1}$:
\begin{equation}
\label{eq:JellettII}
\frac{m}{2}{\area}_{\II}(M)=\int_M H_{\II}\,\rho \,\dd\Omega_{\II}\,.
\end{equation}
\end{theorem}
\begin{proof}[of Theorem~\ref{thm:integral}]
By deforming $M$ in the direction of the position
vector field, we obtain the hypersurface
\[
M_{t}=\left\{n+t P_{(n)}\,\vert\, n\in M \right\}\,\qquad \textrm{(with $t\in \,\left]\,-1\,,\,+\infty\,\right[\,$),}
\]
which is a rescaling of $M$ with a factor $1+t$. Such a homothety magnifies the Gauss-Kronecker curvature $\dA$ with a
factor $\left(\frac{1}{1+t}\right)^m$, whereas $m$-dimensional areas are multiplied with a factor $(1+t)^m$. Hence there holds
$\area_{\II}(M_{t})=(1+t)^{m/2} \area_{\II}(M)$ and consequently,
\[
\deltatnought \area_{\II}(M_{s}) = \frac{m}{2} \area_{\II}(M)\,.
\]
If, on the other hand, $\left\{M_{t}\right\}$ is considered as a deformation of $M$ with
deformation vector field $P$, the definition of the mean curvature of the second fundamental form gives us
\[
\deltatnought \area_{\II}(M_{t})=  \int_M H_{\II}\langle -N, P \rangle
\,\dd\Omega_{\II}\,.
\]
The integral formula is a consequence of these two last equalities.
\end{proof}
\begin{remark}
It is interesting to notice that two integral formulae which appear in  Minkowski's work \cite{minkowski1903}  have been obtained already by J.H.\ Jellett \cite{jellett1853} (see also \cite{bonnet1853}).
The support function plays an important r\^{o}le in these formulae, and the fact that this function has been described already in W.R.\ Hamilton's article \cite{hamilton1833}, suggests that it might have been Hamilton's influence which made the other Dublin mathematician find the formulae. Extensions of these integral formulae can be found in, e.g., \cite{hsiung1954,li1983}, and the integral formula (\ref{eq:JellettII}) may be seen as a special case of the latter one. The integral formula (\ref{eq:JellettII}) can also be obtained by specialising the generalisation of Minkowski's integral formula within the framework of relative differential geometry (see, e.g., \cite{simon1967,suss1929a})
to Manhart's $\II$-normal field.
\end{remark}


\section{Variational Characterisations of Euclidean Spheres.}
\label{sec:varcharac}

In this section, I show that the spheres are the only critical points of the area functional of the second fundamental form under three different constraints.

\subsection{First Variational Problem.}
The proof of our first variational characterisation of the sphere (Theorem~\ref{thm:firstvariationalproblem} below) makes use of the Jellett--Minkowski integral formulae.


\begin{theorem}
\label{thm:firstvariationalproblem}
If an ovaloid in $\mathbb{E}^3$ is a critical point of the
the $\II$-area with respect to variations under which the 
integral mean curvature $\int H\,\mathrm{d}\Omega$ is preserved, then it is a sphere.
\end{theorem}
\begin{proof}[of Theorem~\ref{thm:firstvariationalproblem}]
Under a deformation $\mu$ of an ovaloid $M\subseteq\mathbb{E}^3$ with variational vector field $X$, there holds
\[
\left\{
\begin{array}{rcl}
\displaystyle \deltatnought\, \int_{\mu_t(M)}\! H \,\mathrm{d}\Omega &\displaystyle=&\displaystyle \int_M K \langle -N,X \rangle  \,\mathrm{d}\Omega\,;\\
&&\\
\displaystyle \deltatnought \area_{\II}(\mu_t(M)) &\displaystyle=&\displaystyle \int_M H_{\II}\sqrt{K}\langle -N,X \rangle  \,\mathrm{d}\Omega\,.
\end{array}
\right.
\]
Now assume an ovaloid $M$ is a critical point of the area functional of the second fundamental form, under the constraint that the integral mean curvature be preserved.
The Euler-Lagrange equation which is satisfied on $M$, reads $H_{\II}=C\,\sqrt{K}$ (for $C\in\mathbb{R}$). It should be observed that the constant $C$ is necessarily greater or equal than
one\footnote{This fact can also be deduced from Theorem \ref{thm:ineq}.}:
\[
C\int \sqrt{K}\,\dd\Omega_{\II} = \int H_{\II}\,\dd\Omega_{\II} =
\int H \,\dd\Omega_{\II} \geqslant \int \sqrt{K}\,\dd\Omega_{\II}\,.
\]
On the other hand, an application of formula (\ref{eq:JellettII}) and of another Jellett--Minkowski integral formula gives
\begin{eqnarray*}
\int \sqrt{K}\,\dd\Omega&=&\area_{\II}=  \int H_{\II} \,\rho \,\dd  \Omega_{\II} = C \int K\, \rho  \,\dd  \Omega \\µ
&=&  C\int H \,\dd\Omega \geqslant  \int H \,\dd\Omega \geqslant \int \sqrt{K} \,\dd\Omega \,.
\end{eqnarray*}
This is only possible if $M$ is a sphere.
\end{proof}

\begin{remark} As follows from the inequality
\[
\int H\,\mathrm{d}\Omega \geqslant \int \sqrt{K}\,\mathrm{d}\Omega = \area_{\II}\,,
\]
the spheres are actually a maximum of the functional $\area_{\II}$ which has been constrained to the class consisting of all ovaloids for which the integral mean curvature has some particular value.
\end{remark}

\begin{remark}
As an immediate consequence of Theorem \ref{thm:HIIenH}, the spheres are the only ovaloids in $\mathbb{E}^3$ satisfying
$H_{\II}=C H$ for some constant $C$. The following corollary of Theorem~\ref{thm:firstvariationalproblem} is of a similar nature:
\end{remark}

\begin{corollary}
\label{cor:HIICsqrtK}
The spheres are the only ovaloids in $\mathbb{E}^{3}$ which satisfy $H_{\II}=C\sqrt{K}$ for some $C\in\mathbb{R}$.
\end{corollary}

\begin{remark}
Theorem~\ref{thm:firstvariationalproblem} can be seen as a modification of the following classical theorem: 
``\textit{If an ovaloid in $\mathbb{E}^3$ is a critical point of the (classical) area with respect to variations under which the 
integral mean curvature $\int H\,\mathrm{d}\Omega$ is preserved, then it is a sphere.}''
Namely, the critical points satisfy $H=C K$, and by a theorem of E.B. Christoffel (\cite{christoffel1865}, p. 163) only the spheres satisfy this equation.
We shall give a proof of another modification of Christoffel's result (Theorem \ref{thm:HIICKII} below) which also follows from
theorem 1.d of G.\ Stamou \cite{stamou1987} under the additional assumption $K_{\II}>0$.
\end{remark}

\begin{theorem}
\label{thm:HIICKII}
The spheres are the only ovaloids $M\subseteq \mathbb{E}^3$ such that
\begin{equation}
\label{eq:HIICKII}
\qquad H_{\II}= C K_{\II}
\end{equation}
for some $C\in\mathbb{R}$.
\end{theorem}
\begin{proof}[of Theorem~\ref{thm:HIICKII}]
Let $M$ be an ovaloid such that $H_{\II}=C\,K_{\II}$.
It should be noticed that the constant which occurs in (\ref{eq:HIICKII}) satisfies $C\geqslant 1$. This inequality follows from
\begin{eqnarray}
\label{eq:ch2_HIIisconstKII}
4 \,\pi\, C &=& \int C\, K_{\II} \,\dd\Omega_{\II} = \int H_{\II}\,\dd\Omega_{\II}\\
\nonumber
&=&\int H\,\dd\Omega_{\II} = \int H\,\sqrt{K}\,\dd\Omega \geqslant \int K \,\dd\Omega = 4\,\pi\,.
\end{eqnarray}
The following inequalities are valid at a point $p$ where $K$ achieves its maximal value:
\begin{eqnarray*}
0 &\geqslant& \left\lgroup \frac{1}{4} \Delta_{\II} \log K \right\rgroup_{(p)}
  = \left\lgroup H_{\II}- H \right\rgroup_{(p)} = \left\lgroup C\, K_{\II} - H \right\rgroup_{(p)} \\
  &=& \left\lgroup(C-1)\,H + \frac{C}{2}\mathcal{P} \right\rgroup_{(p)} \geqslant (C-1)\,H(p) \,. 
\end{eqnarray*}
This is only possible if $C=1$, such that equality occurs in (\ref{eq:ch2_HIIisconstKII}). This is only possible if $H=\sqrt{K}$ throughout on $M$.
\end{proof}

\begin{remark}
For an adaption of Corollary~\ref{cor:HIICsqrtK} for curves in the plane, see p.\ 132 of \cite{verpoort2008}. The corresponding result can be reworded as follows: ``\textit{If a simple closed plane curve with strictly positive curvature is a critical point of the length functional of the second fundamental form ($=\int\sqrt{\kappa}\,\dd s$) with respect to deformations under which the (classical) length is preserved, then it is a circle.}''
This leads us towards a second variational problem.
\end{remark} 

\subsection{Second Variational Problem.}

A basic ingredient in the proof of the following theorem is the differential equation 
\[
\Hess_\psi = \frac{\Delta\psi}{m} \,g
\]
for a function $\psi\in\mathfrak{F}(M)$ on an $m$-dimensional Riemannian manifold $(M,g)$. 

The above differential equation naturally arises in the context of conformal transformations between Einstein spaces, between space forms, and concircular transformations (see \cite[lemma 2, proposition 3 and proposition 8, respectively]{kuhnel1988}) and was studied initially by H.W.\ Brinkmann, A.\ Fialkow, Y.\ Tashiro, K.\ Yano \textit{et al}; we give reference to W.\ K\"uhnel's survey text \cite{kuhnel1988}.

\begin{theorem}
\label{thm:secondvariationalproblem}
If an ovaloid in $\mathbb{E}^3$ is a critical point of the 
$\II$-area with respect to deformations under which the (classical) area is preserved, then it is a sphere.
\end{theorem}
\begin{proof}[of Theorem~\ref{thm:secondvariationalproblem}]
An ovaloid is a critical point of this variational problem if and only if the following relation is satisfied for some $C\in\mathbb{R}$\,:
\begin{equation}
\label{eq:ch2_HIIsqrtKH}
H_{\II} \sqrt{K}= C H  \,.
\end{equation} 

Let us temporarily define an \textit{abstract ovaloid} 
as a compact two-dimensional Riemannian manifold  $(M,g)$ such that $M$ is diffeomorphic to a sphere and $g$ has strictly positive Gaussian curvature.
The notation
$\mathcal{E}_{\textsc{abs}} = \left\{\,\textrm{abstract ovaloids}\,\right\}$ will be adopted.

It should be remarked that an abstract ovaloid is actually not so abstract, since the Weyl embedding problem has been solved by H.\ Weyl, H.\ Lewy, \textit{et al}. More precisely, for every abstract ovaloid $(M_{\textsc{abs}},g)$ an ovaloid $M\subseteq\mathbb{E}^3$ can be found which is isometric with $(M_{\textsc{abs}},g)$. Moreover, the congruence theorem of S.\ Cohn-Vossen gives us that $M$ is unique up to a congruence of $\mathbb{E}^3$.

For any deformation $g^{(t)}$ of the metric of an abstract ovaloid $(M,g)$, as well $(M,g^{(t)})$ is an abstract ovaloid for sufficiently small $|t|$, and hence variational problems can be posed on the class $\mathcal{E}_{\textsc{abs}}$. Because of the previous remark, this is essentially not different from the study of variational problems on the class of ``\textit{real}'' ovaloids. 

We introduce the following functionals:
\[
\left\{
\begin{array}{ccccl}
\displaystyle
\mathcal{F} &:& \mathcal{E}_{\textsc{abs}}\rightarrow \mathbb{R} &:& (M,g)\mapsto \displaystyle\int_M \dd\Omega\,;\\
&&&&\\
\displaystyle
\mathcal{F}_{\II} &:& \mathcal{E}_{\textsc{abs}}\rightarrow \mathbb{R} &:& (M,g)\mapsto \displaystyle\int_M \sqrt{K}\,\dd\Omega\,.\\
\end{array}
\right.
\]
Of course, $K$ stands for the Gaussian curvature of the metric $g$ on $M$, and $\dd\Omega$ is the area element of this abstract metric.

If  a deformation $g^{(t)}=g+t\,h+\mathcal{O}(t^2)$ of the metric of an abstract ovaloid $(M,g)$ has been given, then there exists a unique operator $\sigma$, which is symmetric with respect to the metric $g$, and which satisfies $h(V,W) = g(\sigma(V),W)$ for all $V,W\in\mathfrak{X}(M)$.
The following variational formulae can be verified (see, e.g., \cite[pp. 66 ff.]{verpoort2008}, for some details):
\[
\!\!
\left\{\!\!
\begin{array}{rcl}
\displaystyle
\deltatnought\!\!\! \F (M,g^{(t)})
&\!\!\!\displaystyle=&\!\!\! \displaystyle\frac{1}{2} \int_M\! \mathrm{trace}\,\sigma\,\dd\Omega\,;\\
&&\\
\displaystyle\deltatnought\!\!\!\F_{\II} (M,g^{(t)})
 &\!\!\!\displaystyle=&\!\!\!\displaystyle \frac{1}{4} \int_M\!\! 
\mathrm{trace}\left\{\!\left\lgroup 
\left(\sqrt{K}-\Delta\left(\frac{1}{\sqrt{K}}\right)\right)\mathrm{id}
+
\mathrm{Hs}_{\left(\frac{1}{\sqrt{K}}\right)}
\right\rgroup\!\circ\sigma\right\}\,\dd\Omega\,.\\
\end{array}
\right.
\]

By an application of the lemma which is stated on p.\ 168 of \cite{blair2000}, it can be concluded that an abstract ovaloid is a critical point for $\F_{\II}$ under the constraint that $\F$ be preserved if and only if the following operator-valued equation is satisfied:
\begin{equation}
\label{eq:ch_var_ELeq}
\left(\sqrt{K}-\Delta\left(\frac{1}{\sqrt{K}}\right)\right)\mathrm{id}
+
\mathrm{Hs}_{\left(\frac{1}{\sqrt{K}}\right)}
=
C\,\mathrm{id}
\qquad\qquad
\textrm{($C\in\mathbb{R}$).}
\end{equation}

Let now a ``\textit{real}'' and non-spherical ovaloid $M\subseteq\mathbb{E}^3$ be given. Assume $M$ is a critical point of the $\II$-area with respect to deformations under which the (classical) area be preserved. Then the abstract ovaloid $(M,\I)$ is a critical point of $\F_{\II}$ under the constraint that $\F$ be preserved, and consequently (\ref{eq:ch_var_ELeq}) is satisfied.

(If the trace of the composition of both sides of the above equation with the shape operator is taken, there results (\ref{eq:ch2_HIIsqrtKH}), which particularly shows that the constants $C$ occurring in (\ref{eq:ch2_HIIsqrtKH}) and (\ref{eq:ch_var_ELeq}) are equal.)

If the trace of both sides of the above equation (\ref{eq:ch_var_ELeq}) is taken, there results
\begin{equation}
\label{eq:psi}
2\left(\frac{1}{\psi}-C\right) = \Delta \psi\,,
\end{equation}
where $\psi$ stands for the function $\psi=\frac{1}{\sqrt{K}}\in\mathfrak{F}(M)$. 

It follows from the above equation (\ref{eq:ch_var_ELeq}) that
\begin{equation}
\label{eq:ch_var_trace_ELeqbis}
\textrm{Hess}_{\psi} = \frac{\Delta \psi}{2}\, g\,.
\end{equation}
Let $p_+$ be a point where $K$ attains its global maximum (and hence $\psi$ its global minimum). Similarly the global minimum of $K$ (and hence the global maximum of $\psi$)
is attained in a point $p_-$ of $M$. The notation
\[
k_+ = \left.\sqrt{K}\right\vert_{p_+} = \left.\frac{1}{\psi}\right\vert_{p_+}
\qquad\textrm{and}\qquad
k_- = \left.\sqrt{K}\right\vert_{p_-} = \left.\frac{1}{\psi}\right\vert_{p_-}
\]
will be adopted. In the remainder of this proof it will be shown that the assumption that $M$ is non-spherical, which is equivalent to the inequality
\begin{equation}
\label{eq:ch_3_nonspherical}
 k_- < k_+\,,
\end{equation}
allows a contradiction to be deduced.

According to theorem 21 of \cite{kuhnel1988}, $p_+$ and $p_-$ are the only critical points of $\psi$ on $M$. Moreover, according to lemma 22 of \cite{kuhnel1988}, $M\setminus \left\{p_-,p_+\right\}$ (with the first fundamental form metric) is isometric to a warped product 
$\left]\,u_0\,,\,u_1\,\right[\,\times \mathrm{S}^1(\ell)$, the slices of which correspond to the level sets of the function $\psi$. Thus the function $\psi$ can be seen as a function of the first factor only, \textit{i.e.}, $\psi : \,\left]\,u_0\,,\,u_1\,\right[ \,\rightarrow \mathbb{R}$, and, according to  \cite{kuhnel1988}, the derivative of $\psi$ is exactly the warping function. Thus if $\mathrm{S}^1(\ell)$ is described with a co-ordinate $v\in\, \left]\,0\,,\,2\,\pi\,\ell\,\right[\,$, there holds 
\[
M\setminus \left\{p_-,p_+\right\} 
=
\,\left]\,u_0\,,\,u_1\,\right[\,\times_{\psi'}\, 
\left]\,0\,,\,2\,\pi\,\ell\,\right[\,,
\]
with the end-points of the second interval identified.

Then $M\setminus \left\{p_-,p_+\right\}$ can be described with co-ordinates $(u,v)$ and the (first fundamental form) metric is, more precisely,
\begin{equation}
\label{eq:metric}
g = 
\left\lgroup
\begin{array}{cc}
g_{u\,u} & g_{u\,v}\\
g_{u\,v} & g_{v\,v} 
\end{array}
\right\rgroup
=
\left\lgroup
\begin{array}{cc}
1 & 0\\
0 & \left(\psi'(u)\right)^2
\end{array}
\right\rgroup
=
\dd u^2 + \left(\psi'(u)\right)^2 \,\dd v^2\,.
\end{equation}
After possibly having reflected the $u$-interval, the function $\psi : 
\left]\,u_0\,,\,u_1\,\right[\,\rightarrow\mathbb{R}$ is strictly increasing, and there holds
\[
\psi(u_0):=
\lim_{u\rightarrow u_0} \psi(u) = \left.\frac{1}{\sqrt{K}}\right\vert_{p_+} \!\!= \frac{1}{k_+}
\quad
\textrm{and}
\quad
\psi(u_1):=\lim_{u\rightarrow u_1} \psi(u) = \left.\frac{1}{\sqrt{K}}\right\vert_{p_-}\!\!= \frac{1}{k_-}
\,.
\]

The expression (\ref{eq:metric}) for the metric implies that $\displaystyle \frac{\Delta \psi}{2} = \psi''$. Together with (\ref{eq:psi}) this yields the following ordinary differential equation for $\psi$:
\[
\psi'' =\frac{1}{\psi} - C\,.
\]
In particular, since $\lim_{u\rightarrow u_0} \psi(u)$ is a well-defined number in $\mathbb{R}_0^+$, also 
\begin{equation}
\label{eq:chap_var_ode_psi}
\lim_{u\rightarrow u_0} \psi''(u)
=
k_+ - C
\end{equation}
is a well-defined real number (which, of course, is denoted simply by $\psi''(u_0)$).

Now three global geometric invariants can be calculated with respect to this explicit description of the metric on $M$. First of all, the total Gaussian curvature can be expressed as
\begin{equation}
\label{eq:chap_var_intKdA}
4\pi = \int K\,\dd\Omega = \int_{u_0}^{u_1} \int_{0}^{2\pi\ell} \frac{\psi'(u)}{\left(\psi(u)\right)^2}\,\dd v\,\dd u
= 2\,\pi\,\ell \left(k_+ - k_-\right)\,.
\end{equation}
Secondly, the $\II$-area of $M$ satisfies
\begin{equation}
\label{eq:chap_var_areaII}
\area_{\II}(M) = \int \sqrt{K}\,\dd\Omega = 
\int_{u_0}^{u_1} \int_{0}^{2\pi\ell} \frac{\psi'(u)}{\psi(u)}\,\dd v\,\dd u
=
2\,\pi\,\ell \, \ln \left(\frac{k_+}{k_-}\right)\,.
\end{equation}
Thirdly, the area of $M$ satisfies
\begin{equation}
\label{eq:chap_var_area}
\area(M) = \int \dd\Omega = 
\int_{u_0}^{u_1} \int_{0}^{2\pi\ell} \psi'(u)\,\dd v\,\dd u
=2\,\pi\,\ell\,\frac{\left(k_+ - k_-\right)}{k_+\,k_-} \,.
\end{equation}

As immediately follows from a combination of (\ref{eq:ch2_HIIsqrtKH}) together with two Jellett--Minkowski integral equations, these last two quantities are related by 
\begin{equation}
\label{eq:chap_var_areaII_area}
\area_{\II}(M) = \int H_{\II}\, \rho\,\dd\Omega_{\II}
= \int H_{\II}\,\sqrt{K} \,\rho\,\dd\Omega
= C \int H \,\rho\,\dd\Omega
= C \, \area(M)\,.
\end{equation}

Now it will be shown that the constant $C$ satisfies 
\begin{equation}
\label{eq:chap_var_Ck+k-}
C = \frac{\left(k_+ + k_- \right)}{2}\,.
\end{equation}
This can be concluded by a calculation of the length of the geodesic circle $\gamma$ of radius $\varepsilon$ and centre $p_+$ in two different ways.
Firstly, the curve $\gamma$ is described in the co-ordinate system $(u,v)$ by the equation $u=u_0+\varepsilon$, and it can be seen that
\[
\length (\gamma) = \int_{0}^{2\pi\ell} \!\!\psi'(u_0+\varepsilon)\,\dd v = 
2\,\pi\,\ell\,\psi'(u_0+\varepsilon) = 2\,\pi\,\ell\,\psi''(u_0)\,\varepsilon+\Order(\varepsilon^2)\,.
\] 
Secondly, since $\gamma$ is a geodesic circle of radius $\varepsilon$,
\[
\length (\gamma) =
2\,\pi\,\varepsilon + \Order(\varepsilon^2)\,.
\]
A comparison of the two last expressions gives that 
\begin{equation}
\label{eq:chap_var_psidotdotell}
\psi''(u_0)\,\ell=1\,.
\end{equation}
Consequently, making use respectively of (\ref{eq:chap_var_intKdA}), (\ref{eq:chap_var_psidotdotell}) and (\ref{eq:chap_var_ode_psi}), it can be seen that
\begin{eqnarray*}
\ell\,\left(k_+ - k_-\right) 
= 2 
= 2\,\psi''(u_0)\, \ell
= 2 \, \left( k_+ - C \right)\,\ell\,,
\end{eqnarray*}
which gives the promised (\ref{eq:chap_var_Ck+k-}).

By combining equations 
(\ref{eq:chap_var_areaII}),
(\ref{eq:chap_var_areaII_area}),
(\ref{eq:chap_var_area})
and
(\ref{eq:chap_var_Ck+k-}) respectively, we obtain
\begin{equation}
\label{chap_var_equationinkpluskminus}
\ln\left(\frac{k_+}{k_-}\right)
=
\frac{\area_{\II}(M)}{2\,\pi\,\ell}
=
\frac{C\,\area(M)}{2\,\pi\,\ell}
=
\frac{\left(k_+ + k_-\right)\left(k_+ - k_-\right)}{2\, k_+\, k_-}\,.
\end{equation}

Now let the function $f$ be given by 
\[
f: \mathbb{R}_0^+ \rightarrow \mathbb{R}
: x \mapsto \ln(x)-\frac{x}{2}+\frac{1}{2x}\,.
\]
A consideration of the derivative of $f$ suffices to conclude that $f(x)=0$ if and only if $x=1$ (for $x\in\mathbb{R}_0^+$). It has been shown in (\ref{chap_var_equationinkpluskminus}) that $f(\frac{k_+}{k_-})=0$, whence the conclusion $k_+ = k_-$ can be drawn. But this contradicts (\ref{eq:ch_3_nonspherical}), and the theorem is proved.
\end{proof}

\begin{remark}
The spheres are actually a maximum of the constrained functional, as follows from
\[
4\pi\,\area
=
\int K\,\dd\Omega\,\int\,\dd\Omega 
\geqslant
\left\lgroup \int\sqrt{K}\,\dd\Omega \right\rgroup^2
=
\left\lgroup\area_{\II} \right\rgroup^2
\,.
\]
\end{remark}

\begin{remark}
We can also state the following corollary of Theorem~\ref{thm:secondvariationalproblem}:
\end{remark}

\begin{corollary}
The spheres are the only ovaloids in $\mathbb{E}^{3}$ which satisfy $H_{\II}\sqrt{K}=C\,H$ for some $C\in\mathbb{R}$.
\end{corollary}

\subsection{Third Variational Problem.}

The following theorem, which has been mentioned already by M.\ Wiehe in \cite[remark 5.6.(iii)]{wiehe1998}, solves the so-called ``$\II$-isoperimetric problem'' which was formulated in, a.o., \cite{glassner,glassnersimon}. 
The proof of the theorem will be presented as well because of its restricted length: it consists of a mere combination of the work of K.\ Leichtweiss and F.\ Manhart.

\begin{theorem}[K.\ Leichtweiss, F.\ Manhart, M.\ Wiehe]
\label{thm:thirdvariationalproblem}
If an ovaloid in $\mathbb{E}^3$ is a critical point of the 
$\II$-area with respect to deformations under which the volume is preserved, then it is a sphere.
\end{theorem}
\begin{proof}[of Theorem~\ref{thm:thirdvariationalproblem}]
Critical points of this variational problem need to satisfy the relation
$H_{\II} \sqrt{K}= C$ (for some $C\in\mathbb{R}$), or yet
\[
\frac{1}{2}\,\mathrm{trace}\, A_{(N_{\II})} = C\,.
\]

On the other hand, by an application of \cite[Satz 6.1]{simon1967} (see also \cite{suss1929a}), it follows that the relative shape operator with respect to Manhart's $\II$-normal field is a constant multiple of the identity:
\[
\frac{-1}{C}\,\DD_V N_{\II}  = V
\]
for all $V\in\mathfrak{X}(M)$. Since also $\DD_V P = V$, it follows that $\frac{-1}{C}\,N_{\II}= P$ with respect to a suitably chosen origin, and in particular, 
\begin{equation}
\label{eq:leichtweiss}
\sqrt{K}=C\rho\,.
\end{equation}
Since the spheres are the only ovaloids which satisfy the equation 
(\ref{eq:leichtweiss}), as has been shown in theorem 1 of \cite{leichtweiss1995}, the proof is finished.
\end{proof}

\begin{remark}
The spheres are actually a maximum of the constrained functional, as follows from the ``$\II$-isoperimetric inequality'' \cite{manhart1989b}:
\[
48\pi^2\,\vol
\geqslant
\left\lgroup
\area_{\II}
\right\rgroup^3
\,.
\]
\end{remark}

\begin{remark}
We can also state the following corollary of Theorem~\ref{thm:secondvariationalproblem}:
\end{remark}

\begin{corollary}
The spheres are the only ovaloids in $\mathbb{E}^{3}$ which satisfy $H_{\II}\sqrt{K}=C$ for some $C\in\mathbb{R}$.
\end{corollary}


\textbf{Acknowledgement.}
It is a pleasure to thank Dr.\ S.\ Haesen and Professor U.\ Simon for  useful suggestions.


\end{document}